\documentclass[12pt,a4]{article}
\usepackage{amssymb,amsmath}

  \newcommand{\nat}{\mathbb{N}}

\begin{document}
\title{rings in which power values of $K$-Engels with derivations annihilate a certain element }

\author{shervin sahebi$^1$ and venus rahmani$^{2}$\\
Department of Mathematics, Islamic Azad University,\\
Tehran Centre, Tehran, IRAN
\\
e-mail$^1$: sahebi@iauctb.ac.ir\\
e-mail$^2$:ven.rahmani.math@iauctb.ac.ir}

\date{}
\maketitle
\begin{abstract}
let $R$ be a 2-torsion free semiprime ring and $d$  a  non-zero derivation.
Further let $A=O(R)$ be the orthogonal completion of $R$ and
$B=B(C)$ the Boolean ring of $C$ where $C$ be
the extended centroid of $R$.
We show that if  $a[[d(x),x]_{n},[y,d(y)]_{m}]^{t}=0$ such that $0\neq a\in R$
 for all $x,y\in R$, where $m,n,t>0$ are fixed integers,
then there exists an idempotent $e\in B$ such that $eA$
 is a commutative ring and $d$ induce a zero derivation on $(1-e)A$.
\end{abstract}

\noindent
{\bf Math. Subj. Classification 2010:} 16R50; 16N60; 16D60.
\medskip

\bigskip

\baselineskip=20pt

\begin{center} {\bf 1. Introduction}
\end{center}
Let $R$ be an associative ring with center $Z(R)$.
Recall that an additive mapping $d$ of $R$ into itself is a derivation
if $d(xy) = d(x)y + xd(y)$, for all $x, y \in R$.
Also if $(x_{i})_{i\in \nat}$ is a squence of elements of $R$ and $k$
is a positive integer, we define $[x_{1},\ldots,x_{k+1}]$
inductively as follows:
$$[x_{1},x_{2}] = x_{1}x_{2}-x_{2}x_{1}\qquad,\qquad [x_{1},\ldots,x_{k},x_{k+1}] = [[x_{1},\ldots,x_{k}],x_{k+1}].$$
If $x_{1}=x$ and $x_{2}=\ldots=x_{k+1}=y$,
the notation $[x,y]_{k}$ is used to denote $[x_{1},\ldots,x_{k+1}]$
and $[x,y]_{k}$ is called a $k$-Engel element.

\noindent
A well known result  of Posner stated that
if $[[d(x),x],y]=0$ for all $x,y \in R$, then $R$ is commutative~\cite{a8}.
A number of authors  extended this result  in several ways.
Bell and Martindale in~\cite{a0} studied this identity for
a semiprime ring $R$. They proved that if $R$ is a semiprime ring
and $[[d(x),x],y]=0$ for all $x$ in a non-zero
left ideal of $R$ and $y\in R$, then $R$ contains a non-zero central ideal.
In ~\cite{a4}, Filippis showed that if $R$ is a
prime ring with $char R \neq 2$ and $d$ a non-zero
derivation of $R$ such that $[[d(x),x],[d(y),y]]=0$
 for all $x,y\in R$, then $R$ is commutative.
Recently Dhara obtained results
for a prime ring $R$ of $char R\neq 2$,
with a nonzero derivation $d$ that if
$0\neq a\in R$ such that  $a[[d(x),x]_n,[d(y),y]_m]=0$
for all $x, y \in R$, where $m,n\geq 0$ are fixed integers,
 then $R$ is commutative~\cite{a22}.

\noindent
 Now, we will generalize  Posner's result\cite{a8}
 when the condition are more widespread.

  \noindent
 The main result of this paper is as follows:
 
 \noindent
{\bf Theorem 1.1.}\label{THM}
{\it let $R$ be a $2$-torsion free semiprime ring with non-zero
derivation $d$ and $0\neq a \in R$ such that
$a[[d(x),x]_{n},[y,d(y)]_{m}]^{t}=0$ for all $x, y \in R$,
 where $m, n, t> 0$ are fixed integers.
 Further let $A=O(R)$
be the orthogonal completion of $R$ and $B=B(C)$
 where $C$ the extended centroid of $R$.
Then there exists an idempotent $e\in B$ such that $eA$
 is a commutative ring and $d$ induce a zero derivation on $(1-e)A$.}

\noindent
 Throughout the paper we use the standard notation
 from~\cite{a1}. In particular, we denote by  $Q$
 the two sided Martindale quotient of prime and semiprime ring $R$
 and $C$ the center of $Q$. We call $C$ the extended centroid
 of $R$.
 \noindent
It is well known that any derivation of prime(semiprime) ring
$R$  can be uniquely extended to a derivation of $Q$, and so any
derivation of $R$ can be defined on the whole of $Q$.
 Moreover $Q$ is a prime(semiprime) ring as well as $R$.
 We refer to~\cite{a1, a666} for more details.
\bigskip
\begin{center} {\bf 2.Proof of main result} \end{center}
 The following results are usefull tool needed the
 proof of main result.
{\bf Theorem 2.1.}\label{b1}
{\it Let $R$ be a prime ring of $char R\neq 2$ and $d$ a derivation of $R$.
Suppose  $a[[d(x),x]_{n},[d(y),y]_{m}]^{t}=0$ and $0\neq a \in R$
for all $x,y\in R$, where $m, n, t > 0$ are fixed integers.
 Then $R$ is commutative or $d=0$.}

\noindent
{\it Proof.}
 Consider  two cases.
 
\noindent
\emph{case 1}. $d$ is not a $Q$-inner derivation.
 By Kharchenko's Theorem~\cite{a5}
for any $x, y, z, s\in R$ we have
$a[[z,x]_{n},[s,y]_{m}]^{t}=0$.
This is a polynomial identity and hence there
 exists a field $F$ such that $R\subseteq M_{k}(F)$
 with $k>1$ and $R, M_{k}(F)$
satisfy the same polynomial identity ~\cite{a6}.
Therefore we can consider $a=(a_{ij})_{k\times k}$.
We may assume that $t$ is an even integer.
 Now putting
$ z=e_{ij}$, $x=e_{ii}$, $s=e_{ji}$, $y=e_{ii}$.
Thus for any $i\neq j$, we have
$$0=a[[z,x]_{n},[s,y]_{m}]^{t}=a(-1)^{nt}(e_{ii}+(-1)^te_{jj})=a(e_{ii}+e_{jj}),$$
This implies $a_{ij}=0$ for any $i, j$ $(i\neq j)$,
which is contradiction.

\noindent
\emph{case 2}. $d$ is a $Q$-inner derivation.
 So there exists an element $b\in Q$ such
 that $d(x)=[b,x]$ for all $x\in R.$
Since by~\cite{a2} $Q$ and $R$ satisfy the same
generalized polynomial identities $(GPI)$, hence for any $x, y \in Q$ we have
$a[[b,x]_{n+1},[y,[b,y]]_{m}]^t=0$.
Also since $Q$ remains prime by the primeness of $R$,
replacing $R$ by $Q$ we may assume that $b\in R$
and the extended centroid of
$R$ is just the center of $R$.
Note that $R$ is a centrally closed prime
$C$-algebra in the present situation~\cite{a3}.
If $R$ is commutative, we have nothing to prove.
So, let $R$ be noncommutative.
 Therefore $R$ satisfies a nontrivial $(GPI)$.
 Since $R$ is a centrally closed prime $C$-algebra,
 by Martindale's Theorem~\cite{a7},  $R$ is a strongly primitive ring.
Let $_RV$ be a faitful irreducible left
$R$-module with commuting ring $D=End(_RV)$. By the Density Theorem,
$R$ acts densely on $V_{D}.$
For any given $v\in V$ we claim that $v$ and $bv$ are $D$-dependent.
Assume first that $av\neq 0$. Suppose on the contrary that $v$
and $bv$ are $D$-independent.

\noindent 
If $b^{2}v\in span\{v, bv\}$,
then $b^{2}v = v\alpha + bv\beta$ for some $\alpha, \beta \in D.$ By density
of $R$ in $End(V_{D})$ there exist two elements $x$ and $y$ in $R$ such that
$xv = v$, $xbv = 0$ and  $yv = 0$, $ybv = v$.
Then
$$0 = a[[b,x]_{n + 1},[y,[b,y]]_{m}]^{t}v=(-2)^{mt}av.$$
If $b^2v\notin span\{v, bv\},$ then $\{v, bv, b^2v\}$ are all $D$-independent.
Then  by Density of $R$ in $End(V_D)$ there
exist two elements $x$ and $y$ in $R$ such that
$xv=v$, $xbv=0$, $xb^2v=0$ and $yv=0$, $ybv=0$, $yb^2v=0$.
Therefore  we have
$$0 = a[[b,x]_{n + 1},[y,[b,y]]_{m}]^{t}v=(-2)^{mt}av.$$
Since $char R\neq 2$ we get $av=0,$ a contradiction.
Thus $v$ and $bv$ are D-dependent as claimed.
Assume next that $av=0$. Since $a\neq 0$, we have $aw\neq 0$ for some $w\in V.$
Then $a(v+w)=aw\neq 0.$ Applying the first situation we have
$bw=w\alpha$ and $b(v+w)=(v+w)\beta$,
for some $\alpha, \beta\in D.$ But $v$ and $w$ are clearly $D$-independent,
and so there exist two elements $x$ and $y$ in $R$ such that
$xw=w$, $xv=0$ and $yw=v$, $yv=0$.
Then
$$0 = a[[b,x]_{n + 1},[y,[b,y]]_{m}]^{t} = (-1)^{t(n + 1)} 2^{mt} a (\beta -\alpha)^{2t} w,$$
which implies $\alpha = \beta$ and hence $bv = v\alpha$ as claimed.
From the above we have proved that $b v = v \alpha (v)$
for all $v\in V$, where $\alpha (v) \in D$ depends on $v\in V$.
In fact, it is easy to check that $\alpha(v)$ is independent
of the choice of $v\in V.$ That is,
there exist $\delta \in D$ such that $bv = v \delta$ for all $v \in V.$
we claim $\delta \in Z(D)$, the center of $D$. Indeed, if $\beta \in D$, then
$b(v\beta)=(v\beta), \delta=v(\beta\delta)$
and the other hand
$b(v\beta)=(bv)\beta=(v\delta)\beta=v(\delta\beta)$.
Therefore $v(\beta \delta - \delta \beta)=0$ so $\beta \delta=\delta \beta$,
which implies $\delta \in Z(D)$. Thus $b\in C$
and hence $d=0$, as be wanted.

\noindent
 The following example shows the hypothesis of primeness is essential in Theorem $2.1$.
 
 \noindent
{\bf example 2.2.}
{\it Let $S$ be any ring, and
$R =\left\{{\small{\small{\small\left (
      \begin{array}{c c c}
        0 & a & b \\
        0 & 0 & c \\
        0 & 0 & 0
        \end{array} \right )}}}|a, b, c \in S \right\}.$
Define $d:R\rightarrow R$ as follows:
$d {\small{\small\left (
        \begin{array}{c c c}
        0 & a & b \\
        0 & 0 & c \\
        0 & 0 & 0
        \end{array} \right )}}=
 {\small{\small{\small\left (
        \begin{array}{c c c}
        0 & 0 & b \\
        0 & 0 & 0 \\
        0 & 0 & 0
        \end{array} \right )}}}.$
Then $d$ is a non-zero derivation of $R$ such that
$a[[d(x),x]_{n},[d(y),y]_{m}]^{t} = 0$ for
all $x,y \in R$, where $m, n, t > 0$ are fixed integers,
however $R$ is  not commutative.}

Now let $R$ be a  semiprime orthogonally complete
ring with extended centeroid $C$.
We use the notation $B=B(C)$ and spec$(B)$ to
denote Boolian ring of $C$ and the set of all maximal ideal of $B$.
It is well known that if $M\in spec(B)$ then $R_{M}=R/RM$ is prime
~\cite[Theorem 3.2.7]{a1}.
We refer to ~\cite[ pages 37, 38, 43, 120]{a1}
for definations of $\Omega$-$\Delta$-ring,
 a first order formula of signature $\Omega$-$\Delta$,
 Horn formulas and Hereditary first order formulas.

\noindent
In preparation for the proof of Theorem \ref{THM} we have the following lemma.

\noindent
{\bf lemma 2.3.}\label{b2}\cite[Theorem 3.2.18]{a1}.
{\it Let $R$ be an orthogonally
complete $\Omega$-$\Delta$-ring with extended centroid $C$,
$\Psi_{i} ( x_{1}, x_{2},..., x_{n})$ Horn formulas
of signature $\Omega$-$\Delta$, $i=1,2,...$ and
$\Phi(y_{1}, y_{2},...,y_{m})$
a Hereditary first order formula such that
$\neg\Phi$ is a Horn formula. Further, let
$\vec{a}= (a_{1}, a_{2},...,a_{n})\in R^{(n)},$
$\vec{c} = (c_{1}, c_{2},..., c_{m})\in R^{(m)}. $
Suppose that $R\models \Phi (\vec{c})$ and  for every $M\in spect (B)$
there exists a natural number $i=i(M)>0$ such that
$$R_{M} \models \Phi (\phi_{M} (\vec{c})) \Longrightarrow \Psi_{i}(\phi_{M}(\vec{a})),$$
where $\Phi_{M}: R\rightarrow R_{M} = R/RM$
is the canonical projection.
Then there exist a natural number $k > 0$
and pairwise orthogonal idempotents $e_{1}, e_{2},...,e_{k}\in B$ such that
$e_{1} + e_{2} + ... + e_{k} = 1$ and $e_{i}R\models\Psi_{i}(e_{i}\vec{a})$
for all $e_{i}\neq 0$.}

\noindent
Denote by $O(R)$ the orthogonal completion of $R$ which is
defined as the intersection of all orthogonally complete subset
of $Q$  containing $R$.
\smallskip
\noindent
 Now we can prove Theorem 1.1.

\vspace{3mm}
\emph{Proof of Theorem 1.1.}
It is well known that the derivation $d$
can be extended uniquely to a derivation $d: Q\rightarrow Q$.
According to \cite[Theorem 3.1.16]{a1}
$d(A)\subseteq A$ and $d(e)=0$ for all $e\in B.$ Therefore $A$ is an orthogonally
complete $\Omega$-$\Delta$-ring where $\Omega= \{o, +, -, \cdot, d \}$. Consider formulas
$$\begin{array}{l}
\Phi = (\exists a \neq 0) (\forall x ) (\forall y)  \| a[[d(x),x]_{n},[y,d(y)]_{m}]^{t}=0 \|,\\ \\
\Psi_{1} = (\forall x ) (\forall y)  \| xy = yx \|,\\ \\
\Psi_{2} = (\forall x )  \| d(x) = 0 \|.
\end{array}$$
One can easily check that $\Phi $ is a
 hereditary first order formula and $\neg \Phi$,
 $\Psi_{1}$, $\Psi_{2}$ are Horn formulas.
So using Theorem $2.1$ shows that all
conditions of Lemma $2.3$ are fulfilled.
Hence there exist two orthogonal
 idempotent $e_{1}$ and $e_{2}$ such that $e_{1} + e_{2} = 1$ and if $e_{i} \neq 0$,
 then $e_{i}A \models \Psi_{i}$, $i = 1, 2$. The proof is complete.\hfill $\Box$

\end{document}